\documentclass[12pt]{article}

\usepackage{natbib}
\usepackage{amsmath,amssymb}
\usepackage{graphicx}
\usepackage{caption}
\usepackage{subcaption}

\newcommand{\ph}{\ensuremath{\varphi}}
\newcommand{\lap}{\ensuremath{\tilde{\triangle}}}
\newcommand{\ep}{\ensuremath{\varepsilon}}

\title{The Kadomtsev–Petviashvili equation in conformal variables for waves over topography}

\author{David Andrade, \and Marcelo V. Flamarion}

\begin{document}
\maketitle

\begin{abstract}
  The conformal mapping approach is a well established technique for solving the Euler equations for potential flows with one spatial dimension. In this work, we extend this framework to problems with a weakly transversal dependence and, by means of asymptotic expansions, obtain a Kadomtsev–Petviashvili type equation formulated in conformal variables as a model for weakly transversal surface waves propagating over topography. A key advantage of this formulation is that the topography, defined in the physical domain, does not need to be a smooth function, or even a function in the classical sense because, our asymptotic analysis relies on the effective depth, which comes through the Jacobian of the conformal map which is assumed to be a slowly varying function.  The resulting equation provides a consistent extension of several well known weakly nonlinear dispersive wave models previously reported in the literature. Numerical simulations are performed to illustrate the newly derived equation.

\end{abstract}



\section{Introduction}
The conformal mapping method for solving two dimensional potential flow problems has proven to be a powerful analytical and numerical tool in hydrodynamics. It has been successfully applied to a variety of problems, including wave propagation over rough topography \citep{Nachbin2003}, tsunami generation \citep{Poletto2025}, turbulence \citep{Kochurin2025JETP}, the emergence of coherent structures in deep water \citep{Dutykh2013}, flows over bathymetries \citep{Flamarion2019, Ruban2025},  among many others. A comprehensive review about the current state of the art is given in  \citep{Nachbin2023}. 

A principal limitation of the method is that it is naturally formulated for flows with a single horizontal spatial variable. In \citep{AN2018}, the authors propose an extension of the conformal mapping framework to a fully two dimensional free surface flow, on a domain that is independent of the second horizontal variable. Despite significant advances in computational power, fully three dimensional numerical simulations of free surface flows remain computationally demanding. Consequently, reduced models remain valuable tools, as they can capture the essential physical mechanisms while providing tractable descriptions of the dynamics within appropriate asymptotic regimes.

The goal of this work is to exploit the mathematical formulation introduced in \citep{AN2018} and, through an asymptotic expansion, derive a Kadomtsev–Petviashvili type equation (KP equation) in conformal variables. We show that the resulting equation constitutes a natural extension of several well known models in the literature for weakly nonlinear dispersive waves. Two different KP equations will be obtained; one valid for small amplitude topographies and a more sophisticated one for slowly varying depths. Our asymptotic analysis for a small amplitude topography uses a recent expansion proposed by \cite{ludu} whereas the one used for a slowly varying depth appears to be new.. Within this framework, numerical simulations are performed to investigate the interaction of localized pulses and to assess the influence of bathymetric variations on their propagation.

\section{Governing equations}
\label{sec:headings}

Our starting point is the fully non-linear Potential Theory equations for an ideal, three-dimensional, incompressible and irrotational flow under the action of gravity. We assume that the fluid domain is bounded below by a ridge-like bathymetry, i.e. a surface that only depends on $x$, representing a submerged trench or mountain, and above by a free surface $z = \eta(x,y,t)$.

Using $l$, a typical wavelength as length scale, $\sqrt{gh}$ the speed of long waves propagating over $h$, as a reference speed and $a$ a typical wave amplitude, the dimensionless form of the equations is the following:
\begin{eqnarray}
\ph_{xx}+\ph_{yy}+\ph_{zz} &=& 0,\ \ \ \ \text{in $-\mu(1+H(x))<z<\ep\mu\eta(x,y,t)$,}\label{Euler:Eq01}\\
\ph_z +\mu\ph_{x}H_x&=&0,\ \ \ \ \text{on $z=-\mu(1+H(x))$,}\label{Euler:Eq02}\\
\eta_{t} + \ep(\ph_x\eta_x + \ph_y\eta_y) - \tfrac{1}{\mu}\ph_z&=&0,\ \ \ \ \text{on $z = \ep\mu\eta(x,y,t)$,}\label{Euler:Eq03}\\
\ph_t + \eta +\tfrac{\ep}{2}\left(\ph_x^2+\ph_y^2+\ph_z^2\right) &=&0,\ \ \ \ \text{on $z = \ep\mu\eta(x,y,t)$.}\label{Euler:Eq04}
\end{eqnarray}
Two dimensionless parameters arise: $\ep = a/h$, and $\mu = h/l$, the nonlinear and dispersion parameter respectively.

Throughout this paper we will assume that the bathymetry always remains submerged below the free surface but, it will not be assumed that $H$ is necessarily a smooth function or a function at all, see for instance figure \ref{fig1}. In such cases \eqref{Euler:Eq02} is replaced by zero normal derivative of the potential at the seabed.




\subsection{Conformal mapping in three dimensions}

In \citep{AN2018} the conformal mapping technique is extended to the three dimensional domain of equations \eqref{Euler:Eq01} - \eqref{Euler:Eq04}. Such an extension of the conformal map is obtained in the following way. Fix an arbitrary $y$ and let $\Omega_{y}$ be a vertical section of the fluid domain at rest, reflected around $z = 0$:
\begin{equation}\label{Conformal:Eq01}
\Omega_{y} = \{(x,y,z)\in\mathbb{R}^3 \ \mid\ -\mu(1+H(x))<z<\mu(1+H(x))\}.
\end{equation}
Note that $\Omega_y$ is invariant under translations along the $y$-axis. 

Let $F_0$ be a conformal mapping from a uniform strip onto $\Omega_y$:
\begin{equation}\label{Conformal:Eq02}
\begin{aligned}
F_0 &: \{\xi+i\zeta \in\mathbb{C} \mid -\mu < \zeta < \mu\} \longrightarrow \Omega_{y},\\
F_0&(\xi+i\zeta) = x(\xi,\zeta)+iz(\xi,\zeta).
\end{aligned}
\end{equation}
Next, extend $F_0$ to a map $F$ from the uniform three dimensional domain onto the (symmetrized) physical fluid domain:
\begin{equation}\label{Conformal:Eq03}
\begin{aligned}
F &: \{(\xi,y,\zeta) \in \mathbb{C}\times\mathbb{R} \mid -\mu < \zeta < \mu\} \longrightarrow \cup_{y\in\mathbb{R}}\Omega_{y}\\
F&(\xi,y,\zeta) = (x(\xi,\zeta),y,z(\xi,\zeta)).
\end{aligned}
\end{equation}
Note that $F$ `flattens out' the fluid domain and that its coordinate functions are harmonic.

\subsection{Potential theory in harmonic coordinates}

Using the coordinate functions of $F$ as a new \emph{harmonic coordinate system}, equations (\ref{Euler:Eq01})-(\ref{Euler:Eq04}) can be rewritten as, see \citep{AN2018}:
\begin{eqnarray}
\ph_{\xi\xi}+J\ph_{yy}+\ph_{\zeta\zeta} &=& 0,\ \ \ \ \text{in $-\mu<\zeta<\ep\mu N(\xi,y,t)$,}\label{EulerConformal:Eq01}\\
\ph_\zeta &=&0,\ \ \ \ \text{on $\zeta=-\mu$,}\label{EulerConformal:Eq02}\\
JN_{t} + \ep(\ph_\xi N_\xi + J\ph_yN_y) - \tfrac{1}{\mu}\ph_\zeta&=&0,\ \ \ \ \text{on $\zeta = \ep\mu N(\xi,y,t)$,}\label{EulerConformal:Eq03}\\
\ph_t + \eta +\tfrac{\ep}{2J}\left(\ph_\xi^2+J\ph_y^2+\ph_\zeta^2\right) &=&0,\ \ \ \ \text{on $\zeta = \ep\mu N(\xi,y,t)$.}\label{EulerConformal:Eq04}
\end{eqnarray}
In these equations we use the following notation, taken from \citep{Nachbin2003}. Let $J(\xi,\zeta) = z^2_\xi(\xi,\zeta) + z^2_\zeta(\xi,\zeta)$ be the Jacobian of the conformal mapping $F_0$, and let $N$ be a scaled pre-image of the free surface $\ep\mu\eta$, in the sense that:
\begin{equation}\label{EulerConformal:Eq05}
\zeta = \ep\mu N(\xi,y,t) = \operatorname{Im}\left(F^{-1}(x,y,\ep\mu\eta(x,y,t))\right),
\end{equation}
where $\operatorname{Im}$ stands for the imaginary part. 

\subsection{A Boussinesq model over ridge-like bathymetries}

The three dimensional extension of the \emph{Terrain Following Boussinesq System} of \citep{Nachbin2003}, is 
\begin{eqnarray}
    &M\eta_t + \lap{\ph} + \ep\left[\frac{\eta}{M}\ph_\xi\right]_\xi + \ep\left[\frac{\eta}{M}\ph_y\right]_y + \frac{\mu^2}{3}I\ph_{yy} = 0.\label{Bou:eq01}\\
    &\ph_t + \eta + \frac{\ep}{2}\left(\left(\frac{\ph_\xi}{M}\right)^2 + \ph_y^2\right) - \frac{\mu^2}{3}\lap\ph_t = 0. \label{Bou:eq02}
\end{eqnarray}

In equations \eqref{Bou:eq01} and \eqref{Bou:eq02} the effects of the bathymetry are incorporated through the new horizontal variable $\xi$, the two variable coefficients $M(\xi)$, $I(\xi)$ and the second order differential operator $\tilde{\triangle}$, which are given by:
\begin{eqnarray}\label{Bou:coefficients}
    &M^2(\xi) = J(\xi,0) = z_
    \xi^2(\xi,0) + z_\zeta^2(\xi,0).\label{M}\\
    &I(\xi) = M^2_\xi - M_{\xi\xi} \label{I}\\
    &\tilde{\triangle} = \partial_{\xi\xi} + M^2(\xi)\partial_{yy}.
\end{eqnarray}
All the details about the derivation of these equations and the coefficients are given in \citep{AN2018}.

\begin{figure}
    \centering
    \includegraphics[width=1\linewidth]{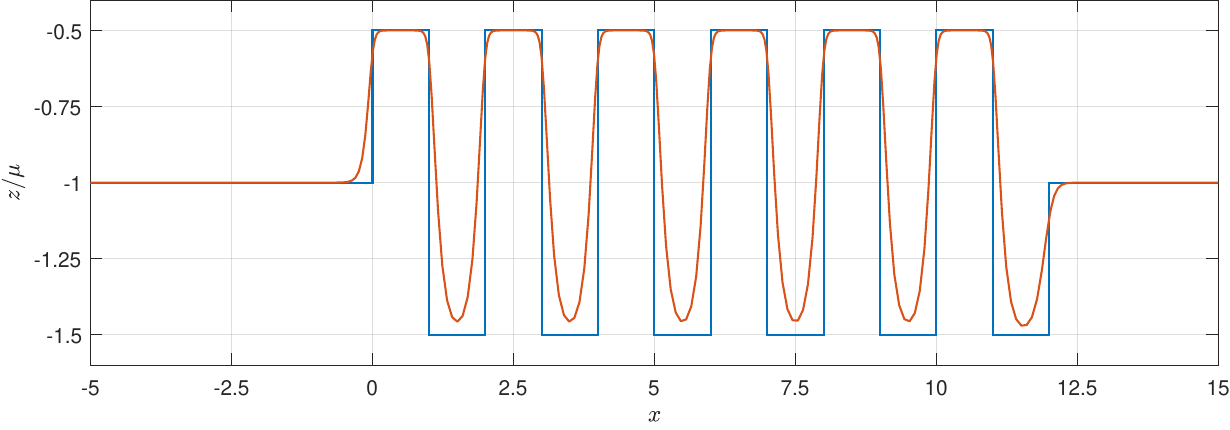}
    \caption{Blue line: Shape of the submerged topography $z/\mu = -(1 + H(x))$. Note that the vertical variable $z$ has been normalized by $\mu = \sqrt{0.1}$. The topography is given by a sequence of rectangles with amplitudes rangind from $-1.5$ to $-0.5$, and all having width 1. Red line: plot of the curve $(x(\xi),-M(\xi))$ which is the variable coefficient that encodes the topographic information of the model.}
    \label{fig1}
\end{figure}

In order to obtain a KP equation for the propagation of weakly nonlinear, weakly dispersive and weakly transversal waves over topography, we must adopt the usual KP scaling whereby the transversal coordinate $y$ is rescaled as $\gamma y$, where $\gamma = \mathcal{O}(\mu)$ and, consistent with the validity of Boussinesq model $\ep = \mathcal{O}(\mu^2)$ the resulting system of equations simplifies to the following:
\begin{eqnarray}
    &M\eta_t + \ph_{\xi\xi} + \gamma^2M^2\phi_{yy} + \ep\left[\frac{\eta}{M}\ph_\xi\right]_\xi= 0.\label{Bou:eq03}\\
    &\ph_t + \eta + \frac{\ep}{2}\left(\left(\frac{\ph_\xi}{M}\right)^2 \right) - \frac{\mu^2}{3}\ph_{\xi\xi t} = 0. \label{Bou:eq04}
\end{eqnarray}
Upon differentiating with respect to $\xi$ and letting $u = \ph_\xi$, the system is recast as
\begin{eqnarray}
    &M\eta_t + u_{\xi} + \gamma^2M^2(\partial^{-1}_\xi u)_{yy} + \ep\left[\frac{\eta}{M}u\right]_\xi= 0.\label{Bou:eq05}\\
    &u_t + \eta_\xi + \frac{\ep}{2}\left(\left(\frac{u}{M}\right)^2\right)_\xi - \frac{\mu^2}{3}u_{\xi\xi t} = 0. \label{Bou:eq06}
\end{eqnarray}

These equations \eqref{Bou:eq05} and \eqref{Bou:eq06} are the starting point of the asymptotic derivation of the KP equation.

\section{A KP equation for slowly varying depth}

 In order to kickstart the asymptotic derivation of the KP equation we assume that the variable coefficient $M(\xi)$ is a slow function of $\xi$. This does not necessarily mean that the topography is a slowly varying function itself but, on the contrary, we should expect that $M(\xi)$ is a slowly varying function whenever the topography is. In general, we can always expect $M(\xi)$ to be a smoother version of the topography profile and that it varies on a slower scale. This follows from the fact that, as shown in \citep{Nachbin2003}, $M(\xi)$ is obtained as a convolution between the topography with a smoother kernel. Our modeling hypothesis is that $M$ depends on the slow conformal variable $\mu^2\xi$ so:

\begin{align}\label{eq:M slow}
    M(\xi) = M(\mu^2\xi).
\end{align}

The derivation of a KP-like equation from equations \eqref{Bou:eq05} and \eqref{Bou:eq06} begins by considering the first-order system
\begin{eqnarray}
    &M\eta_t + u_{\xi} = 0,\label{Bou:eq07}\\
    &u_t + \eta_\xi  = 0. \label{Bou:eq08}
\end{eqnarray}

The Riemann invariants associated with the system \eqref{Bou:eq07}--\eqref{Bou:eq08} satisfy
\begin{equation}
(R_+)_t+\frac{1}{\sqrt{M}}(R_+)_\xi=0
\end{equation}
and
\begin{equation}
(R_-)_t-\frac{1}{\sqrt{M}}(R_-)_\xi=0.
\end{equation}

The corresponding Riemann invariants are given by
\begin{equation}
R_+ = u + \sqrt{M}\,\eta
\end{equation}
and
\begin{equation}
R_- = u - \sqrt{M}\,\eta.
\end{equation}

These quantities remain constant along the characteristic curves
\begin{equation}
\frac{d\xi}{dt}=\frac{1}{\sqrt{M}}
\end{equation}
and
\begin{equation}
\frac{d\xi}{dt}=-\frac{1}{\sqrt{M}},
\end{equation}
respectively.

Since we are interested in right-propagating waves, we impose $R_-=0$, which yields
\begin{equation}
u = \sqrt{M}\,\eta.
\end{equation}

Motivated by this, we seek solutions of the form
\begin{equation} \label{ansatz}
u = \sqrt{M}\eta+\ep A(\xi,y, t) +\mu^2 B(\xi,y,t) +\gamma^2 C(\xi,y,t).
\end{equation}

Replacing \eqref{ansatz} into equations \eqref{Bou:eq05} and \eqref{Bou:eq06}, using our modeling assumption to get the correct scaling of the derivatives of $M$, and collecting the terms of the same order yields:
\begin{eqnarray}
    &M\eta_t + \sqrt{M}\eta_{\xi} +  \ep\Big[A_{\xi} +\frac{2}{\sqrt{M}}\eta\eta_\xi\Big]+\mu^2\Big[B_\xi+\frac{M_\xi}{2\sqrt{M}}\eta  \Big] +\gamma^2\Big[C_\xi +M^2\partial_{\xi}^{-1}(\sqrt{M}\eta_{yy}) \Big]= 0.\label{Bou:eq09}\\
    &\sqrt{M} \eta_t + \eta_\xi +  \ep\Big[A_{t} +\frac{1}{M}\eta\eta_\xi \Big] +\mu^2\Big[B_t-\frac{1}{3}\sqrt{M}\eta_{\xi\xi t}\Big]+\gamma^2 C_t = 0. \label{Bou:eq10}
\end{eqnarray}
Multiply equation \eqref{Bou:eq09} by $\sqrt{M}$ and equate the resulting expression with equation \eqref{Bou:eq10}. This yields the following compatibility conditions
\begin{equation}\label{A}
A_t-\frac{1}{\sqrt{M}}A_\xi=\frac{1}{M}\,\eta\eta_\xi
\end{equation}

\begin{equation}\label{B}
B_t-\frac{1}{\sqrt{M}}B_\xi=
\frac{\sqrt{M}}{3}\eta_{\xi\xi t}
+\frac{M_\xi}{2M}\eta
\end{equation}

\begin{equation}\label{C}
C_t-\frac{1}{\sqrt{M}}C_\xi=
M^{3/2}\partial_\xi^{-1}\!\left(\sqrt{M}\,\eta_{yy}\right)
\end{equation}
Using the relation $\partial_t=-1/\sqrt{M}\partial_{\xi}$, valid to leading order, into the expressions \eqref{A}-\eqref{C} yields
\begin{equation} \label{Ax}
A_\xi=-\frac{1}{2\sqrt{M}}\eta\eta_\xi
\end{equation}

\begin{equation} \label{Bx}
B_\xi
=
-\frac{M_\xi}{4\sqrt{M}}\,\eta
+\frac{1}{6}\sqrt{M}\,\eta_{\xi\xi\xi}
\end{equation}

\begin{equation} \label{Cx}
C_\xi=-\frac{1}{2}M^2\partial_\xi^{-1}\left(\sqrt{M}\eta_{yy}\right)
\end{equation}

Replacing equations \eqref{Ax}–\eqref{Cx} into \eqref{Bou:eq09} yields the Kadomtsev–Petviashvili equation with variable coefficients
\begin{equation}\label{KP1}
\eta_t
+\frac{1}{\sqrt{M}}\eta_\xi
+\frac{3\ep}{2M\sqrt{M}}\eta\eta_\xi
+\frac{\mu^2 M_\xi}{4M\sqrt{M}}\eta
+\frac{\mu^2}{6\sqrt{M}}\eta_{\xi\xi\xi}
+\frac{\gamma^2}{2}M\partial_\xi^{-1}(\sqrt{M}\eta_{yy})
=0.
\end{equation}

Note that, as a consistency check, when the depth is constant, $x = \xi,$ $M = 1$ and the classical Kadomtsev–Petviashvili equation is recovered.

\subsection{Comparison with existing models of wave propagation over slowly varying depth}

In order compare our KP equation \eqref{KP1} with other existing models of wave propagation over slowly varying depth, we must re write our equation in the physical variable $x$.

We can undo the conformal transformation by setting $x = x(\xi,0)$, so $dx/d\xi = M(\xi)$ and the differential operators are $\partial_x = \partial_\xi/M$. Note that the relation $x = x(\xi)$ is invertible and we denote its inverse by $\xi = \xi(x)$ and $d\xi/dx = 1/M(\xi(x))$. We define the \emph{effective depth} as
\begin{equation}
    d(x) = M(\xi(x)),
\end{equation}
so $d\xi/dx = 1/d(x)$.

Setting $\gamma = 0$ in equation \eqref{KP1} yields the following version of the KdV equation 
\begin{equation}\label{KdV1:conformal}
    \eta_t
+\frac{1}{\sqrt{M}}\eta_\xi
+\frac{3\ep}{2M\sqrt{M}}\eta\eta_\xi
+\frac{\mu^2 M_\xi}{4M\sqrt{M}}\eta
+\frac{\mu^2}{6\sqrt{M}}\eta_{\xi\xi\xi}
=0,
\end{equation}
which, in the original physical variables transforms into
\begin{equation}\label{KdV1:physical}
    \eta_t + c\eta_x + \frac{3\ep}{2}\frac{c}{d}\eta\eta_x + \frac{\mu^2}{6}cd^2\eta_{xxx} + \mu^2\frac{c_x}{2}\eta = 0,
\end{equation}
where $c(x) = \sqrt{d(x)}$. This is exactly the KdV equation for a slowly varying depth reported in \citep{Grimshaw1981, Johnson1973a, Johnson1973b} but with the effective depth in place of the physical topography.

The KP equation for slowly varying depth in physical space is
\begin{equation}\label{KP1:physical}
    \partial_x\left(\eta_t + c\eta_x + \frac{3\ep}{2}\frac{c}{d}\eta\eta_x + \frac{\mu^2}{6}cd^2\eta_{xxx} + \mu^2\frac{c_x}{2}\eta\right) + \frac{\gamma^2 c}{2}\eta_{yy} = 0.
\end{equation}

\section{A KP equation for small amplitude depth}

Another KP model for wave propagation that can be obtained from the Boussinesq system \eqref{Bou:eq05} and \eqref{Bou:eq06} occurs when the topography has small amplitude an so we can assume that $M$ is of the form
\begin{align}\label{eq: M small}
    M(\xi) = 1 + \ep m(\xi).
\end{align}
This form of a perturbation expansion in $M$ was first introduced in \cite{ludu}, to obtain a simplified Boussinesq system that allowed the authors to describe asymptotically, the behavior of solitary waves that propagate over a low corrugation periodic topography. We shall use their idea to obtain a KP equation for a small amplitude, ridge-like bathymetry.   

To leading order, we get the first order system
\begin{eqnarray}
    &\eta_t + u_{\xi} = 0,\\
    &u_t + \eta_\xi  = 0, 
\end{eqnarray}
which admits unidirectional waves provided that $u = \eta$. 

We look for solutions of the form 
\begin{align}
    u = \eta + \ep A + \mu^2 B + \gamma^2C.
\end{align}
Then, upon substituting into equations \eqref{Bou:eq05} and \eqref{Bou:eq06}, and equating the terms for $A$ and $B$ respectively, we get the following compatibility conditions
\begin{equation}\label{A:KP small}
A_\xi - A_t = -\eta\eta_\xi - \ep m\eta_\xi.
\end{equation}
\begin{equation}\label{B: KP small}
B_t-B_\xi=
-\frac{1}{3}\eta_{\xi\xi t}.
\end{equation}
\begin{equation}\label{C: KP small}
C_t-C_\xi=
-\partial_\xi^{-1}\left(\eta_{yy}\right).
\end{equation}

Using the fact that, to leading order, $\partial t = -\partial_\xi$ and substituting back into the Boussinesq system we get the following KP equation for small amplitude topography:
\begin{equation}\label{KP2}
    \eta_t + \left(1 - \frac{\ep}{2}m\right)\eta_\xi + \frac{3\ep}{2}\eta\eta_\xi + \frac{\mu^2}{6}\eta_{\xi\xi\xi} + \frac{\gamma^2}{2}\partial_\xi^{-1}\left(\eta_{yy}\right)=0.
\end{equation}

Neglecting the transversal dependence, i.e., setting $\gamma = 0$ in equation \eqref{KP2} yields a KdV equation that is similar to the one obtained  recently in \cite{ludu} on a somewhat different reference frame that ours.

In order to get back to the physical variables a new function called the \emph{effective amplitude} of the topography is introduced as
\begin{equation}
    h(x) = m(\xi(x)).
\end{equation}
Its relation to the effective depth is
\begin{equation}
    d(x) = 1 + \ep h(x).
\end{equation}

Using the same relations between the conformal and physical variables mentioned above, the resulting KP equation for small amplitude topography is
\begin{equation}\label{KP2:physical}
    \partial_x\left(\eta_t + \eta_x + \frac{\ep h}{2}\eta_x + \frac{3\ep}{2}\eta\eta_x + \frac{\mu^2}{6}\eta_{xxx}\right) + \frac{\gamma^2}{2}\eta_{yy} = 0.
\end{equation}

\section{Numerical simulations}

 In this section we illustrate the differences between the KP equations, i.e. equations \eqref{KP1} and \eqref{KP2} by considering the propagation of an incoming wave that propagates over a patch rectangles, exactly as shown in Figure \ref{fig1}. The unidirectional nature of the KP equations makes it suitable for accurate numerical computations by means of pseudo spectral methods, so we impose periodic boundary conditions.

In the following simulations we used a square domain of the form $[-L_x,L_x]\times[-L_y,L_y]$, with $L_x = L_y = 30$. We set up a uniform grid with $1024$ points along each direction, which yields a resolution of $dx = dy = 60/1024\approx 0.058$. In all simulation we set $\ep = 0.05$ and, to be consisten with the model, $\mu = \gamma = \sqrt{0.05}$. As it is customary with spectral methods, differentiations are carried out though the FFT and, in order to avoid aliasing, multiplications are computed back in the physical space. 

\begin{figure}
     \centering
     \begin{subfigure}[b]{1\linewidth}
         \centering
         \includegraphics[width=\linewidth]{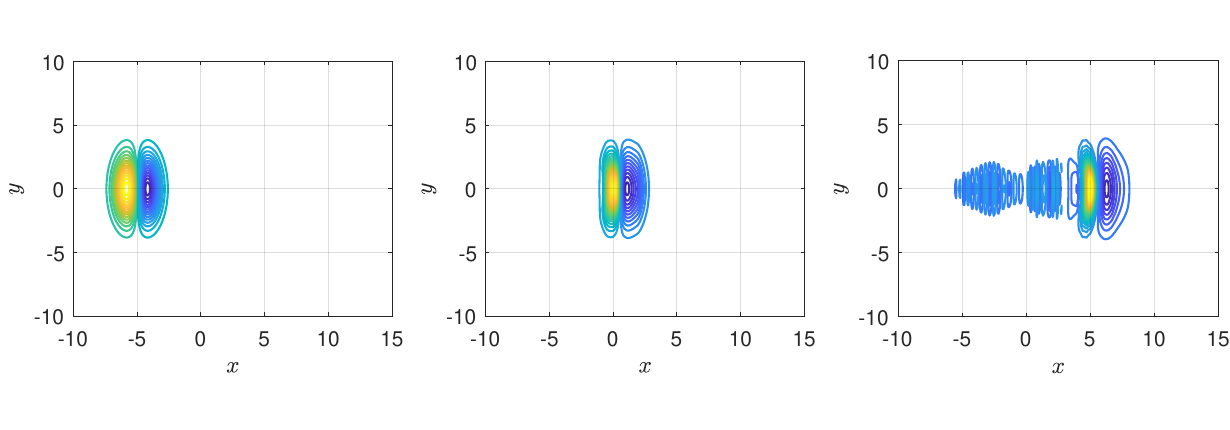}
         \caption{Solution of the KP equation \eqref{KP1} valid for slowly varying $M$.}
         \label{fig2a}
     \end{subfigure}
     \hfill
     \begin{subfigure}[b]{1\linewidth}
         \centering
         \includegraphics[width=\linewidth]{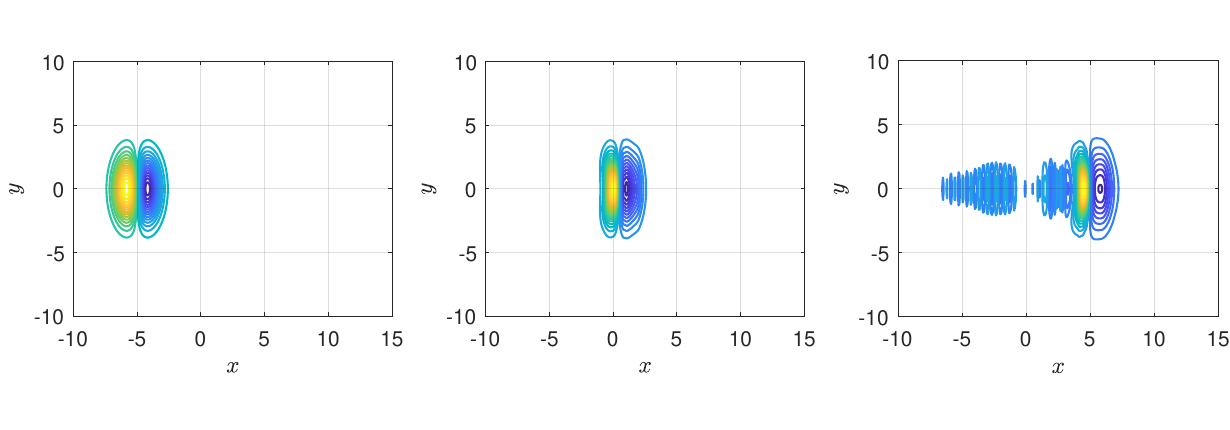}
         \caption{Solution of the KP equation \eqref{KP2} valid for a small amplitude $M$. }
         \label{fig2b}
     \end{subfigure}
     \hfill
     \begin{subfigure}[b]{1\linewidth}
         \centering
         \includegraphics[width=\linewidth]{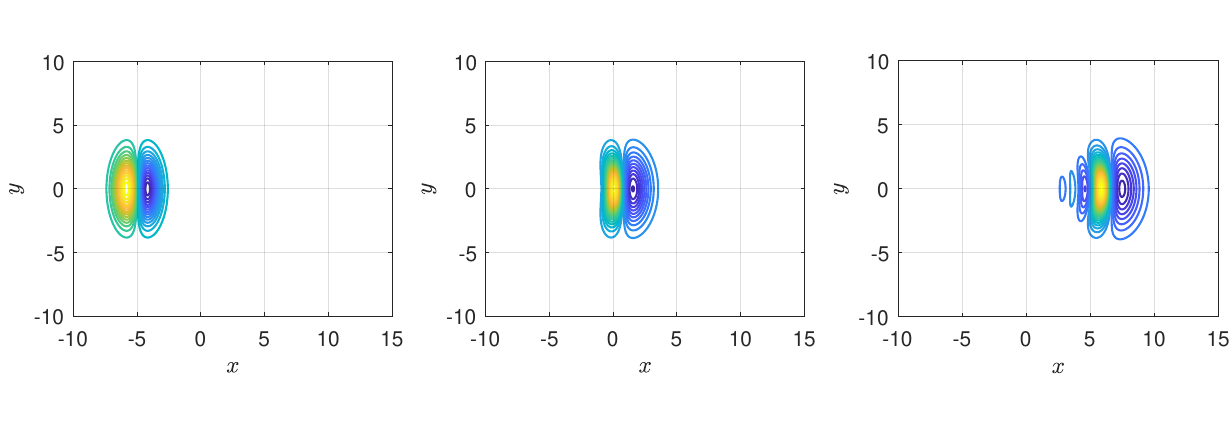}
         \caption{Solution of the classical KP equation, valid for $x = \xi$ and $M = 1$.}
         \label{fig2c}
     \end{subfigure}
        \caption{Snapshots of the time evolution of an incoming wave at times $t = 0$, $6$ and $12$. The topography is located on the interval $0\leq x\leq 12$ and it is the one given in Figure \ref{fig1}.}
        \label{fig2}
\end{figure}

The initial conditions are give by 
\begin{align}
    \eta(x,y,0) = \frac{\partial}{\partial x} \exp\left(-\frac{(x + 5)^2}{2s^2_x}\right)\exp\left(-\frac{y^2}{2s^2_y}\right),
\end{align}
so it is the derivative of a Gaussian pulse with $s_x = 0.82$ and $s_y= 1.65$. Level lines of the initial data are shown in the left panel of Figure \ref{fig2}. Time evolution is computed using MATLAB ode45 with absolute and relative tolerances set to $10^{-12}$.

In Figure \ref{fig2} we show the time evolution of an incoming wave, simulated from the KP equation \eqref{KP1} in Figure \ref{fig2a}, from the KP equation \eqref{KP2} in Figure \ref{fig2b}, and, in order to see the effects of the topography on the wave  propagation, the time evolution from the classical KP equation is also shown in Figure \ref{fig2c}.

From Figure \ref{fig2}, one can see the effect of the variable topography on the wave propagation. First of all, the topography slows down the incoming wave as it is clear from the rightmost panels in figures \ref{fig2a}, \ref{fig2} and \ref{fig2c}. Also one can see that the combined effect of the topography with dispersion leads to an oscillatory trailing wave behind the main pulse. Note that wake is clearly different for both equations. see Figures \ref{fig2a} and \ref{fig2b}.

\section{Discussion and conclusions}

We have obtained two versions of the KP equation for slowly varying depth and their corresponding KdV versions in the conformal variables, which admit the use of `rough' topographies in physical space through the Jacobian of the transformation $M(\xi)$. Following the observation in \citep{Nachbin2003}, one should expect $M(\xi) = M(\mu\xi)$ already, so our modeling assumption $M(\xi) = M(\mu^2\xi)$ should not be too stringent and be met in a many cases of interest even though the topography itself $H(x)$ need not meet such condition. A more thorough and complete argument supporting this point has recently been made in \cite{ludu} and it justifies the use of equations \eqref{eq:M slow} and \eqref{eq: M small} as the starting point of our results.

One key point of this paper is that, although the  submerged topography of the physical model need not be a smooth function, or a function at all, through the conformal map we have obtained what we called the \emph{effective depth} i.e., the smooth function $d(x) = M(\xi(x))$, which is the  facto depth for the reduced KdV or KP model written down in the original physical coordinates. This explains the conventional wisdom that one can actually use `any' depth profile in the existing reduced models for wave propagation over depth, without being too worried about meeting the slowly varying (or small amplitude) hypothesis used in the derivation of the model. Moreover, we have found that `any' depth can actually be used provided that one feeds into those models their respective effective depth.

\textbf{Declaration of interests.} The authors report no conflict of interest
    






\bibliographystyle{apalike}
\bibliography{Referencias.bib}


\end{document}